\documentclass[preprint,12pt,authoryear]{elsarticle}

\usepackage{amssymb}
\usepackage{amsmath}
\usepackage{hyperref}
\usepackage{subcaption}
\usepackage{float}
\usepackage{listings}
\usepackage{xcolor}
\begin{document}

\begin{frontmatter}



\title{Two-dimensional Gauss--Jacobi Quadrature for Multiscale Boltzmann Solvers} %% Article title


\author{Shanshan DONG}{}
\author{Lu WANG}\cortext[cor]{Corresponding author: 
	Lu Wang: wanglu@hdu.edu.cn}
\author{Xiangxiang CHEN}{}
\author{Guanqing WANG}{} %% Author name

\affiliation{organization={Hangzhou Dianzi University},%Department and Organization
            addressline={}, 
            city={Hangzhou},
            postcode={310018}, 
            state={Zhejiang},
            country={China}}

%% Abstract
\begin{abstract}
The discretization of velocity space plays a crucial role in the accuracy and efficiency of multiscale Boltzmann solvers. Conventional velocity space discretization methods suffer from uneven node distribution and mismatch issues, limiting the performance of numerical simulations. To address this, a Gaussian quadrature scheme with a parameterized weight function is proposed, combined with a polar coordinate transformation for flexible discretization of velocity space. This method effectively mitigates node mismatch problems encountered in traditional approaches. Numerical results demonstrate that the proposed scheme significantly improves accuracy while reducing computational cost. Under highly rarefied conditions, the proposed method achieves a speed-up of up to 50 times compared to the conventional Newton-Cotes quadrature, offering an efficient tool with broad applicability for numerical simulations of rarefied and multiscale gas flows.
\end{abstract}

\begin{graphicalabstract}

	\centering
	\includegraphics[width=0.45\textwidth]{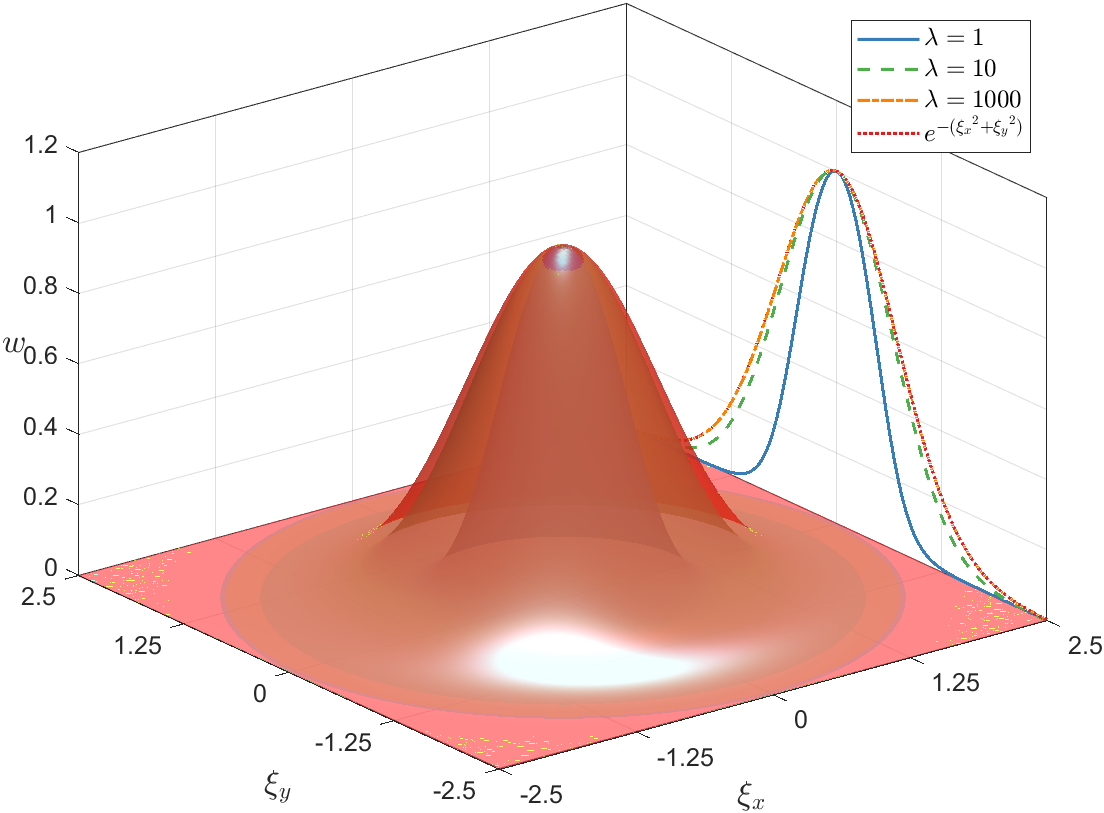}\quad
	\includegraphics[width=0.45\textwidth]{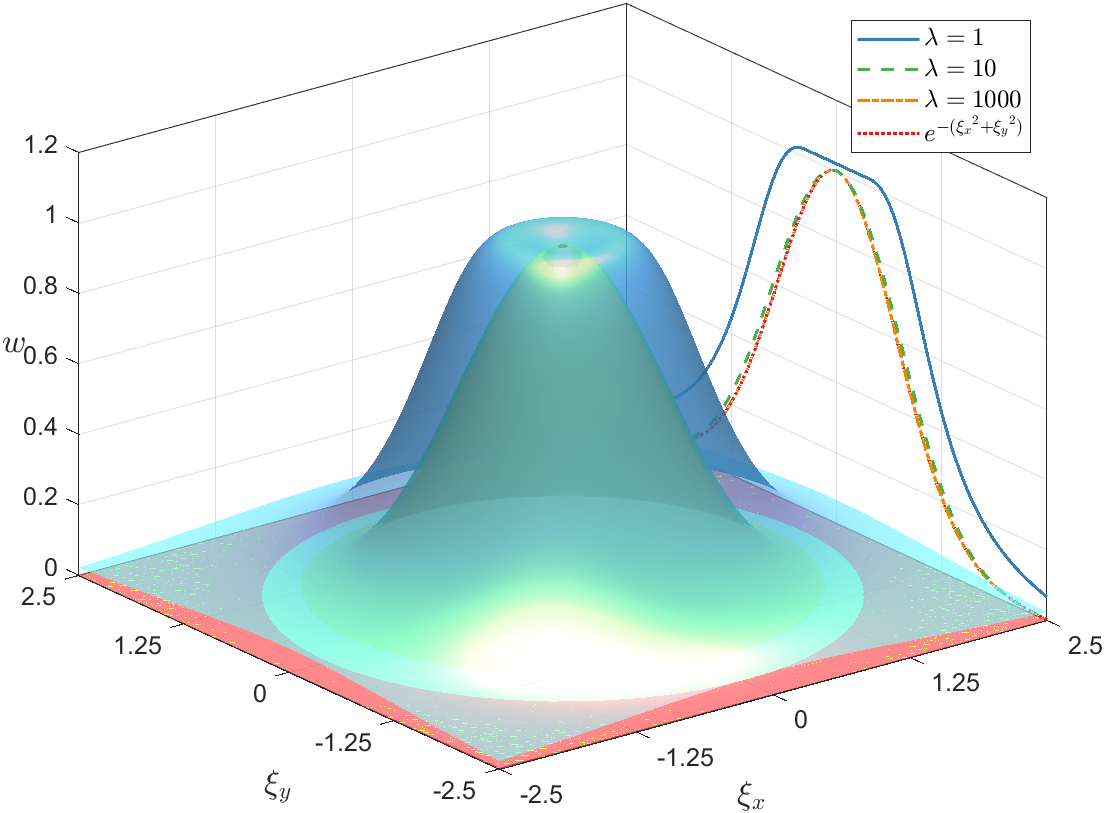}
	
\end{graphicalabstract}

\begin{highlights}

\item Proposed a Gaussian quadrature with parameterized weight and polar transform.

\item Solved node distribution and mismatch issues in velocity space discretization.

\item Achieved higher accuracy with reduced computational cost.

\item Delivered up to $50 \times$ speed-up under highly rarefied conditions.

\end{highlights}

%% Keywords

\end{frontmatter}



\section{Introduction}
\label{sec1}
Accurate simulation of multiscale gas dynamics is increasingly critical in fields such as aerospace engineering, vacuum technology, and micro-electro-mechanical systems (MEMS). To address this challenge, numerical methods based on the discrete velocity method (DVM)\citet{ref1,ref2,ref3} have been extensively developed for solving the Boltzmann equation. Representative schemes include the Gas-Kinetic Unified Algorithm (GKUA)\citet{ref4,ref5}, Unified Gas-Kinetic Scheme (UGKS)\citet{ref6,ref7}, and Discrete Unified Gas-Kinetic Scheme (DUGKS)\citet{ref8,ref9,ref10}, all of which offer unified treatment across the full range of Knudsen numbers. A key component of these methods is the accurate and efficient discretization of the velocity space.

Among various kinetic models, the Shakhov model is widely adopted due to its ability to correct the Prandtl number while preserving energy conservation\citet{ref3,ref6,ref11}.The performance of DVM-based solvers heavily relies on velocity space quadrature. Existing quadrature methods fall into two main categories. The first employs Gaussian-type quadrature over infinite domains, such as Gauss--Hermite and Gauss--Laguerre rules\citet{ref12,ref13}. Although accurate near equilibrium, their fixed node locations limit performance in high-Knudsen or high-Mach-number regimes. The second category truncates the velocity domain and applies finite-interval rules, such as Gauss--Legendre, Gauss--Chebyshev, and Newton--Cotes quadrature\citet{ref13,ref14,ref15}. These methods offer more flexibility but typically require dense grids, resulting in high computational cost.

To improve adaptability, recent efforts have explored parameterized quadrature schemes. Wang et al.\cite{ref16} constructed a novel bell-shaped weight function based on hyperbolic tangent functions and proposed the Gauss--Jacobi quadrature rule on infinite intervals. This method introduces parameters to adjust the discrete velocity distribution and shows good accuracy and adaptability. However, it relies on tensor-product extensions to multidimensional velocity spaces, leading to exponentially increasing computational cost with dimensionality. More recently, a parameteric Gaussian quadrature (PGQ) in polar or spherical coordinates has been introduced to reduce dimensionality and improve efficiency\citet{ref17}. However, the use of a fixed weight function in PGQ still limits node distribution flexibility. Nevertheless, the PGQ method is limited by its fixed weight function, which restricts the flexibility of the discrete velocity distribution and reduces its adaptability in multiscale flow simulations.

In this work, we develop a new Gaussian quadrature with a fully tunable weight function in polar coordinates. The proposed method allows for adaptive control of both node positions and weights, enabling accurate and efficient simulation across a broad range of Knudsen and Mach numbers. Numerical experiments demonstrate improved accuracy and reduced computational cost compared with conventional schemes.
\section{Mathematical Model}
\label{sec2}
This work adopts the BGK--Shakhov model formulation of the Boltzmann equation, which introduces a heat flux correction term to the classical BGK model, thereby enabling a more accurate representation of non-equilibrium heat transfer phenomena under varying Prandtl number conditions. In two-dimensional space, the nondimensional reduced BGK--Shakhov equations are given by~\citet{ref13}:
\begin{eqnarray}
	\frac{\partial g}{\partial t} + \boldsymbol{\xi} \cdot \nabla g = \varOmega(g) \equiv -\frac{1}{\tau}(g - g^S),
\end{eqnarray}
\begin{eqnarray}
	\frac{\partial h}{\partial t} + \boldsymbol{\xi} \cdot \nabla h = \varOmega(h) \equiv -\frac{1}{\tau}(h - h^S),
\end{eqnarray}
where \( g(\boldsymbol{x}, \boldsymbol{\xi}, t) \) and \( h(\boldsymbol{x}, \boldsymbol{\xi}, t) \) are the distribution functions for mass and energy, respectively; \( \boldsymbol{\xi} \) is the particle velocity; \( \tau \) is the relaxation time; and \( g^S \), \( h^S \) are the Shakhov-modified equilibrium distribution functions:
\begin{eqnarray}
	g^S = g^{eq} \left[ 1 + (1 - \text{Pr}) \frac{4 \boldsymbol{c} \cdot \boldsymbol{q}}{5 p T} \left( \frac{|\boldsymbol{c}|^2}{T} - 2 \right) \right],
\end{eqnarray}
\begin{eqnarray}
	h^S = \frac{(1 + N) T}{2} g^{eq} \left[ 1 + (1 - \text{Pr}) \frac{2 \boldsymbol{c} \cdot \boldsymbol{q}}{5 p T} \left( \frac{2|\boldsymbol{c}|^2}{T} - 2 - \frac{2N}{1 + N} \right) \right].
\end{eqnarray}

The equilibrium distribution \( g^{eq} \) takes the form:
\begin{eqnarray}
	g^{eq} = \frac{\rho}{\pi T} \exp\left( -\frac{|\boldsymbol{c}|^2}{T} \right),
\end{eqnarray}
where \( \boldsymbol{c} = \boldsymbol{\xi} - \boldsymbol{u} \) is the peculiar velocity, \( \boldsymbol{u} \) is the macroscopic velocity, \( \boldsymbol{q} \) is the heat flux vector, \( p = \rho T \) is the pressure, \( \rho \) is the density, and \( N \) denotes the internal degrees of freedom.

The macroscopic moments and their basis functions are defined as:
\begin{eqnarray}
	\boldsymbol{M} &= 
	\begin{bmatrix}
		\rho \\
		\rho \boldsymbol{u} \\
		\rho E \\
		\boldsymbol{q}
	\end{bmatrix},
	\quad
	\boldsymbol{\psi} = 
	\begin{bmatrix}
		1 & \boldsymbol{\xi} & \dfrac{1}{2}|\boldsymbol{\xi}|^2 & \dfrac{1}{2} \boldsymbol{c} |\boldsymbol{c}|^2 \\
		0 & 0 & \dfrac{1}{2} & \dfrac{1}{2} \boldsymbol{c}
	\end{bmatrix}^\mathrm{T}.
\end{eqnarray}

The macroscopic variables are then obtained via moment integration:
\begin{eqnarray}
	\boldsymbol{M} = \int \boldsymbol{\psi} 
	\begin{bmatrix}
		g \\
		h
	\end{bmatrix} 
	d\boldsymbol{\xi}.
\end{eqnarray}

The distribution functions $g\left( \boldsymbol{x},\boldsymbol{\xi },t \right)$ and $h\left( \boldsymbol{x},\boldsymbol{\xi },t \right)$ are continuous functions of the particle velocity $\boldsymbol{\xi}$. In the discrete velocity method (DVM), the velocity space is discretized into a set of discrete velocities $\boldsymbol{\xi}_i$, whose specific values are determined by the chosen quadrature rule. Consequently, the quadrature rule determines the number of discrete velocity distribution equations to be solved, and thus directly affects the overall computational cost.
\section{Numerical Quadrature Rules}
\label{sec3}
The velocity distribution function in the Boltzmann equation is defined over an unbounded domain, whereas classical Gauss–Jacobi quadrature is applicable only to finite intervals. To resolve this domain mismatch, we introduce an arctangent transformation that maps the infinite velocity space onto a finite interval, enabling the construction of an arctangent-based Gauss–Jacobi quadrature (ATGJ). The resulting weighted quadrature rule, associated with a parameterized weight function $\omega_{\alpha,\lambda}$, takes the form:
\begin{eqnarray}
	\label{eq:atgj_start}
	I(f) = \int_{\mathbb{R}^2} \omega_{\alpha,\lambda}(\xi_x,\xi_y) f(\xi_x,\xi_y)\, d\xi_x d\xi_y 
	\approx \sum_{k=1}^K w_k f(\xi_{x,k}, \xi_{y,k}),
\end{eqnarray}
where $\xi_{x,k}$ and $\xi_{y,k}$ are the discrete velocity nodes and $w_i$ are the corresponding weights. The weight function is defined as:
\begin{eqnarray}
	\label{eq:omega}
	\omega_{\alpha,\lambda} = \frac{\left[1 - \frac{2}{\pi} \arctan \left( \chi_{\xi} \right) \right]^{\alpha}}{1 + \chi_{\xi}^2},
\end{eqnarray}
with $\chi_{\xi} = \frac{\xi_x^2 + \xi_y^2}{\lambda T_0}$, $\alpha, \lambda > 0$ being tunable parameters.

To derive the ATGJ quadrature, we apply a polar coordinate transformation to Eq.~\eqref{eq:atgj_start}, letting
\begin{eqnarray}
	\label{eq:polar}
	\xi_x = R(r) \cos\theta,\quad
	\xi_y = R(r) \sin\theta,
\end{eqnarray}
where $R(r)=\sqrt{\lambda T_0 \tan\left( \frac{\pi}{2} r \right)}$, $r \in (0,1)$ and $\theta \in (0, 2\pi)$. Substituting Eq.~\eqref{eq:polar} into Eq.~\eqref{eq:atgj_start}, the integral becomes:
\begin{eqnarray}
	\label{eq:polar_transformed}
	I(f) = \frac{\pi}{4} \lambda T_0 \int_0^1 (1 - r)^{\alpha} \int_0^{2\pi} 
	f\left( R(r) \cos\theta, R(r) \sin\theta \right) d\theta\, dr.
\end{eqnarray}

The radial integral over $r$ is evaluated using a Gauss–Jacobi quadrature with weight $(1 - r)^\alpha$, and the angular integral over $\theta$ is computed using a periodic trapezoidal rule:
\begin{eqnarray}
	\label{eq:theta_nodes}
	\theta_j = \theta_0 + \frac{2\pi j}{N_\theta},\quad w_{\theta_j} = \frac{2\pi}{N_\theta},\quad j=1,2,\dots,N_\theta.
\end{eqnarray}

Combining the two, the ATGJ quadrature becomes:
\begin{eqnarray}
	\label{eq:atgj_final}
	I(f) \approx \frac{\pi}{4} \lambda T_0\, w_\theta 
	\sum_{i=1}^{n} w^{\alpha,\lambda}_{r,i}
	\, f\left( R(r_i^{\alpha,\lambda}) \cos \theta_j,\,	R(r_i^{\alpha,\lambda}) \sin \theta_j \right),
\end{eqnarray}
where $r_i^{\alpha,\lambda}$ and $w^{\alpha,\lambda}_{r,i}$ are the Gauss–Jacobi nodes and weights, and $\theta_j$ and $w_{\theta_j}$ are the angular nodes and weights from Eq.~\eqref{eq:theta_nodes}.  
The detailed computation of \(R(r_i^{\alpha,\lambda})\) and \(w^{\alpha,\lambda}_{r,i}\) is provided in the Appendix.  
Accordingly, the discrete velocity nodes and weights for Eq.~\eqref{eq:atgj_start} are:
\begin{eqnarray}
	\label{eq:discrete_nodes}
	\xi_{x,k} = R(r_i^{\alpha,\lambda}) \cos \theta_j,\quad
	\xi_{y,k} = R(r_i^{\alpha,\lambda}) \sin \theta_j,
	\quad
	w_k = \frac{\pi}{4} \lambda T_0\, w^{\alpha,\lambda}_{r,i}\, w_{\theta_j} .
\end{eqnarray}

Although the weight function $\omega_{\alpha,\lambda}$ in Eq.~\eqref{eq:omega} appears complex, it leads to a compact quadrature rule with favorable mathematical properties. Its profile closely resembles a Gaussian distribution and forms a smooth bell-shaped surface. Notably, when $\alpha = \frac{\pi}{2} \lambda$, the weight function asymptotically approaches the classical Maxwellian distribution in the limit $\lambda \to \infty$:
\begin{eqnarray} 
	\label{eq:omega_limit}
	\lim_{\lambda \to \infty} \omega_{\alpha,\lambda} = \exp\left( -\frac{\xi_x^2 + \xi_y^2}{T_0} \right).
\end{eqnarray}

To verify Eq.~\eqref{eq:omega_limit}, let \( k = \xi_x^2 + \xi_y^2 \). As \( \lambda \to \infty \), we have \( \frac{k}{\lambda T_0} \to 0 \). Since the denominator in Eq.~\eqref{eq:omega} tends to unity in this limit, it suffices to analyze the behavior of the numerator. Specifically,
\begin{eqnarray}
	\label{eq:limit_proof}
	\left[ 1 - \frac{2}{\pi} \arctan\left( \frac{k}{\lambda T_0} \right) \right]^{\frac{\pi}{2}\lambda}
	\sim \left( 1 - \frac{2k}{\pi \lambda T_0} \right)^{\frac{\pi}{2} \lambda} 
	\to \exp\left( -\frac{k}{T_0} \right).
\end{eqnarray}

The asymptotic result in Eq.~\eqref{eq:omega_limit} confirms that the classical Maxwellian function is a special case of the proposed weight formulation. This convergence ensures that the ATGJ scheme remains consistent with traditional Gaussian quadrature while offering greater flexibility through its tunable parameters. Such flexibility is particularly advantageous in discretizing velocity space for rarefied and multiscale gas flow simulations.

\begin{figure}[h]
	\centering
	\includegraphics[width=0.45\textwidth]{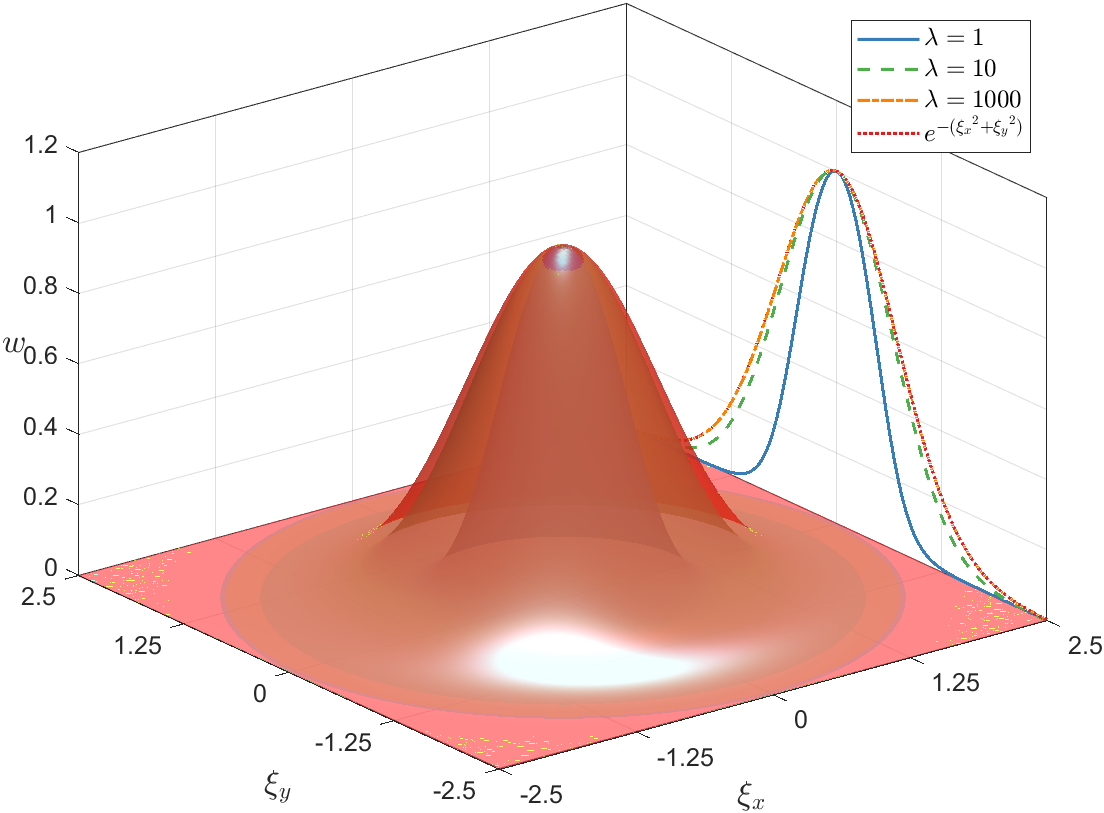}\quad
	\includegraphics[width=0.45\textwidth]{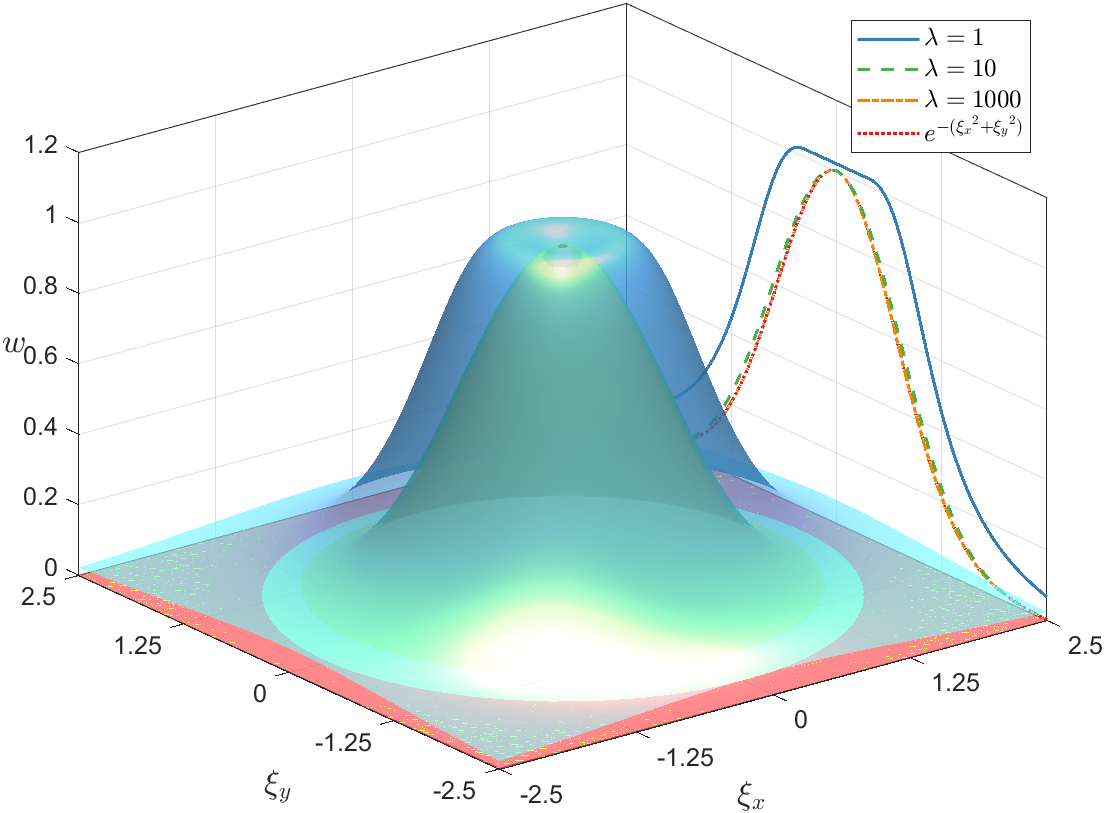}
	
	\vskip 1mm
	
	\makebox[0.45\textwidth][c]{(a)}\quad
	\makebox[0.45\textwidth][c]{(b)}
	
	\caption{
		Variations of bell-shaped weight functions under different parameter settings: 
		(a) $\alpha = \frac{\pi}{2}\lambda + 2$ \quad
		(b) $\alpha = \frac{\pi}{2}\lambda - 2$
	}
	\label{fig:1}
\end{figure}

Fig.~\ref{fig:1} illustrates the surface profile of $\omega_{\alpha,\lambda}$ under various parameter settings with $T_0 = 1$. When $\alpha > \frac{\pi}{2} \lambda$, the bell-shaped surface expands outward and gradually approaches the Maxwellian from within, as shown in Fig.~1(a). Conversely, for $\alpha < \frac{\pi}{2} \lambda$, the surface contracts inward and converges to the Gaussian profile from the outside, as shown in Fig.~1(b). This tunable convergence behavior enables accurate and adaptive coverage of the Maxwellian distribution, supporting efficient and robust quadrature design for kinetic simulations.
\section{Numerical Test Cases}
\label{sec4}
In this section, the effectiveness of the proposed method is validated through two benchmark problems: temperature-discontinuity-induced (TDI) cavity flow under different Knudsen numbers ($\mathrm{Kn}$), and high-Mach-number flow past a square cylinder. These cases span flow regimes from near-continuum to highly rarefied and supersonic conditions. The discrete unified gas kinetic scheme (DUGKS) for the Shakhov model, developed by Guo et al.~\citet{ref8}, is adopted for numerical implementation.

\begin{figure}[htbp]
	\centering
	\begin{minipage}[t]{0.42\textwidth}
		\centering
		\raisebox{-3mm}{\includegraphics[width=\linewidth]{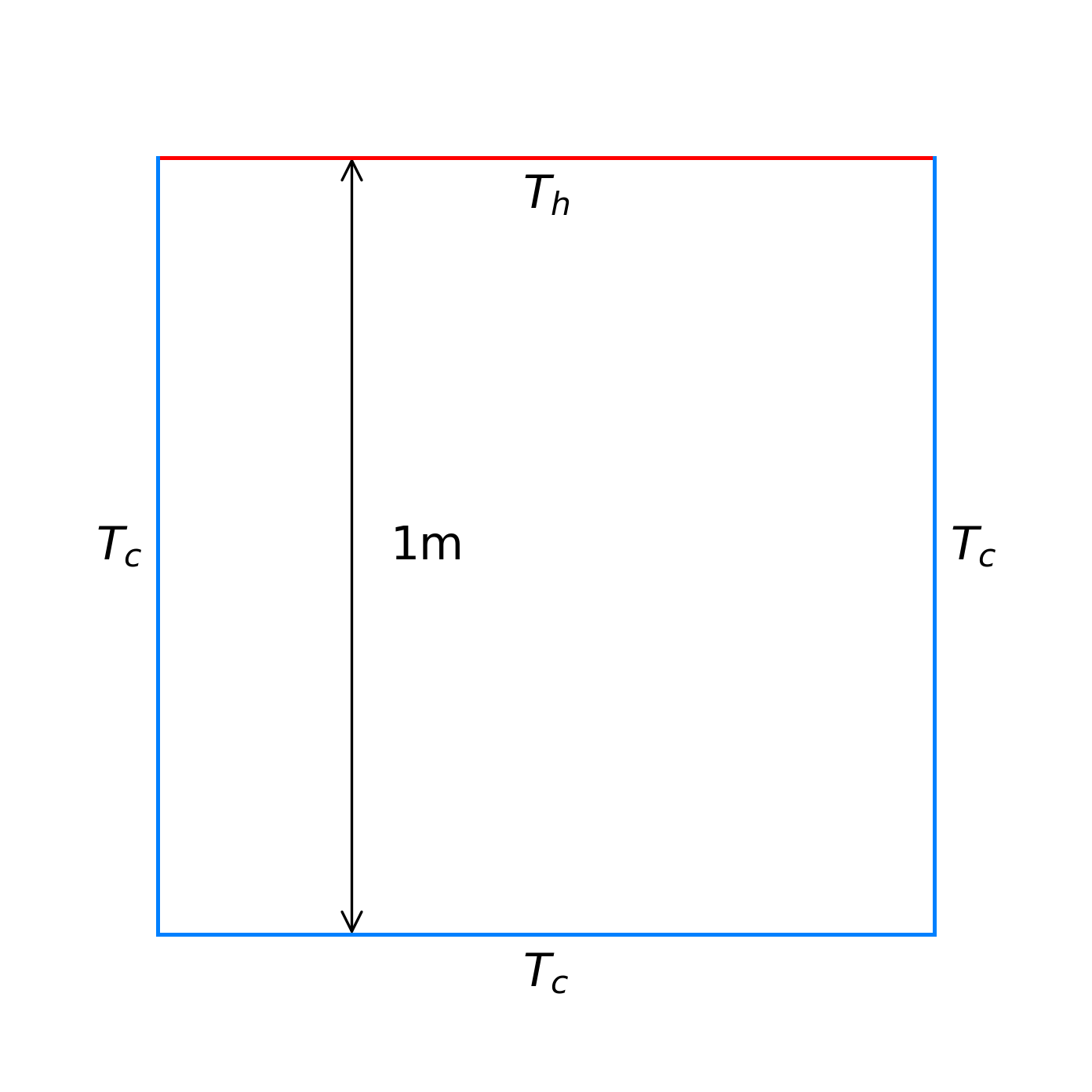}}
		\vskip 0.5mm
		(a)
	\end{minipage}
	\quad
	\begin{minipage}[t]{0.48\textwidth}
		\centering
		\includegraphics[width=\linewidth]{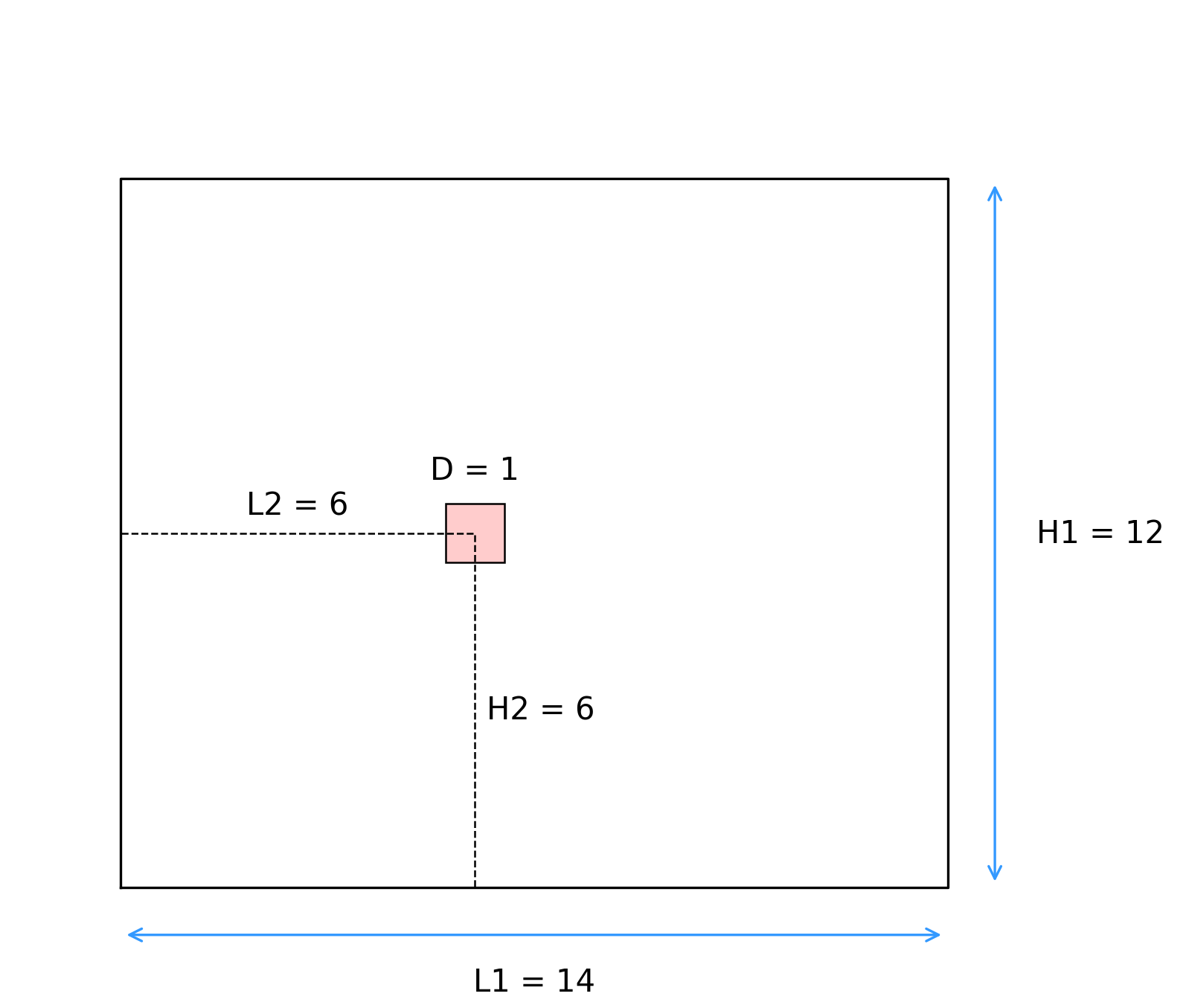}
		\vskip 0.5mm
		(b)
	\end{minipage}
	
	\vskip 2mm
	
	\caption{Illustration of flow geometries: 
		(a) TDI cavity flow; \quad
		(b) Hypersonic flow past a square cylinder.}
	\label{fig:flow_geometries}
\end{figure}

\subsection{TDI cavity flow}
\label{subsec1}
TDI cavity flow is a standard benchmark in multiscale gas dynamics, widely used to investigate natural convection and non-equilibrium heat transport under varying Knudsen numbers. The geometry is shown in Fig.~2(a), where a square cavity of side length $L = 1\,\mathrm{m}$ has its top wall maintained at a higher temperature $T_h$, while the remaining walls are kept at $T_c$. Zhu \textit{et al.}~\citet{ref18} simulated this flow using the DUGKS method across four Knudsen numbers, employing Half-Range Gauss--Hermite (HGH) and Newton--Cotes (NC) quadrature rules. Their velocity discretization settings are summarized in Table~1.
\begin{table}[H]
	\centering
	\caption{Discrete velocity settings for TDI cavity flow under different Knudsen numbers}
	\label{tab:velocity_settings}
	\vskip 2mm
	\footnotesize
	\setlength{\tabcolsep}{4.5pt} 
	\renewcommand{\arraystretch}{1.1} 
	\begin{tabular}{lcccc}
		\hline\hline
		Method & Kn=0.001 & Kn=0.1 & Kn=1 & Kn=10 \\
		\hline
		Zhu et al. & $12 \times 12$ (HGH) & $28 \times 28$ (HGH) & $161 \times 161$ (NC) & $201 \times 201$ (NC) \\
		ATGJ & $8 \times 16$ & $8 \times 45$ & $8 \times 60$ & $8 \times 90$ \\
		Ratio & 1.1 & 2.2 & 54.0 & 56.1 \\
		\hline\hline
	\end{tabular}
	\renewcommand{\arraystretch}{1} 
\end{table}

In the continuum regime ($\mathit{Kn} = 0.0001$), wall temperatures are set to $T_h = 301\,\mathrm{K}$ and $T_c = 300\,\mathrm{K}$. For all other cases, $T_h = 400\,\mathrm{K}$ and $T_c = 200\,\mathrm{K}$ are used. The reference temperature is fixed at $T_{\mathrm{ref}} = 300\,\mathrm{K}$, and fully diffusive boundary conditions are imposed. A uniform $60 \times 60$ mesh is used for spatial discretization. For velocity space, results obtained with the proposed ATGJ scheme are compared to those in Ref.~\cite{ref18}. In the ATGJ formulation, the velocity discretization is controlled via the parameters $\alpha = \frac{\pi}{2} \lambda$. For the near-continuum case, a large value of $\lambda = 500$ is selected to approximate the Gaussian weight, while for rarefied regimes, a smaller value $\lambda = 5$ is used. Under this setting, discrete velocities are distributed approximately within $[-4, 4]$, effectively covering the particle velocity range in low-Mach-number flows.

\begin{figure}[H]
	\centering
	\begin{subfigure}[b]{0.45\textwidth}
		\raisebox{-2mm}{\includegraphics[width=\textwidth]{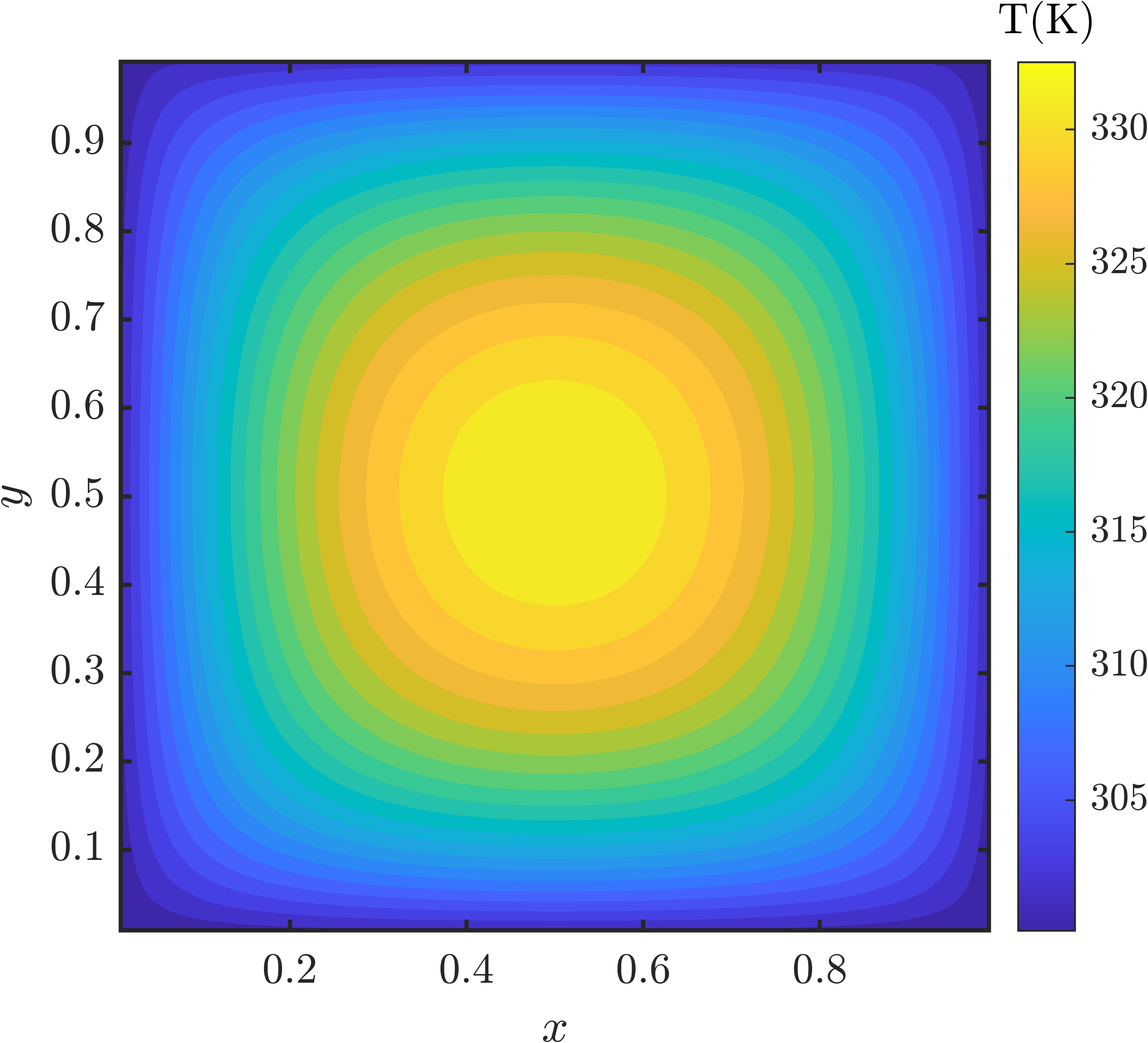}}
		\caption{}
		\label{fig:temp_a}
	\end{subfigure}
	\quad
	\begin{subfigure}[b]{0.45\textwidth}
		\includegraphics[width=\textwidth]{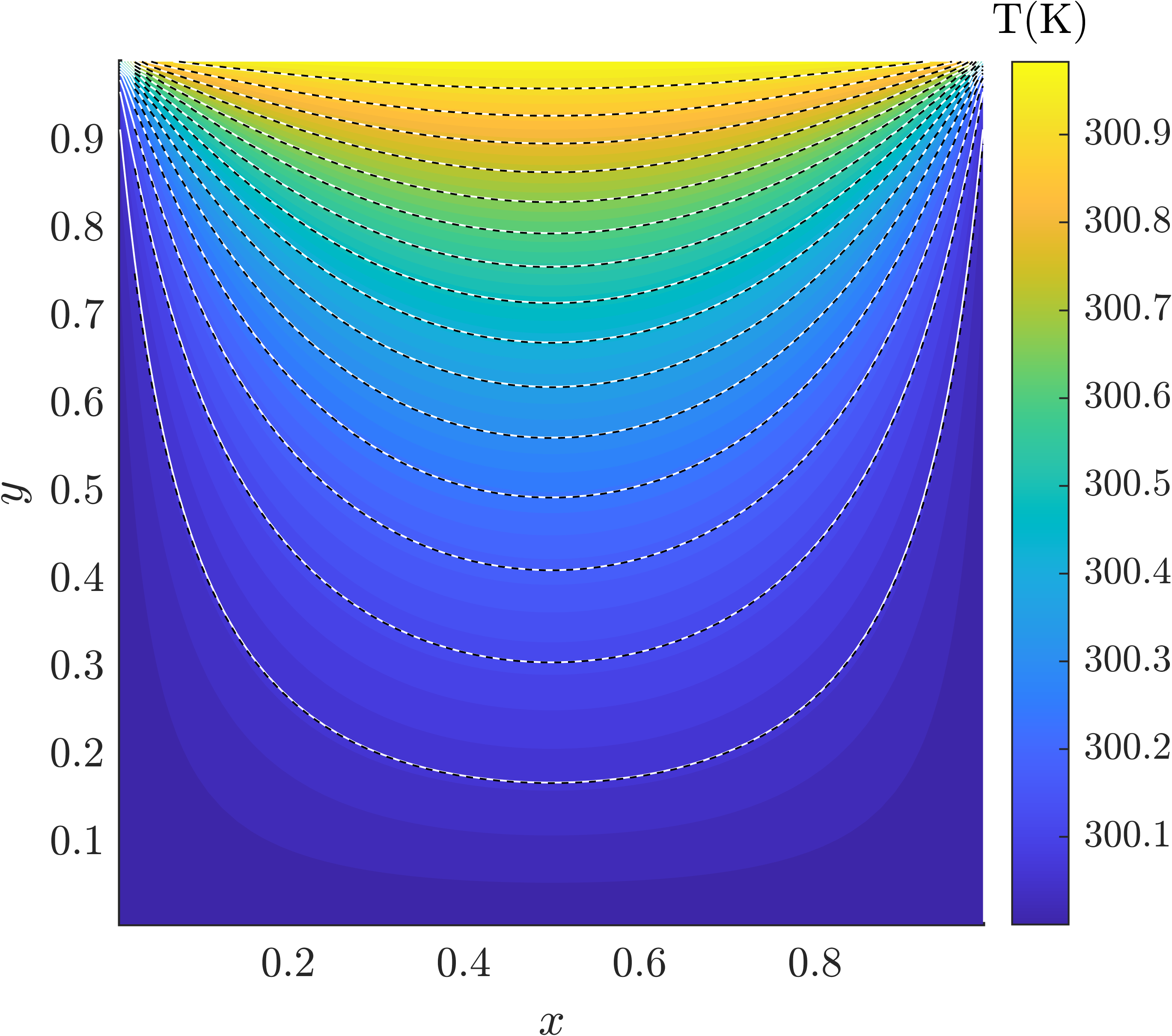}
		\caption{}
		\label{fig:temp_b}
	\end{subfigure}
	\caption{Temperature contours of the TDI cavity flow at $\mathit{Kn} = 0.001$:  
		(a) $\lambda = 5$; (b) $\lambda = 500$. The white solid lines denote the analytical solution, and the black dashed lines correspond to the ATGJ results.}
	\label{fig:temperature_tdi}
\end{figure}

Fig.~\ref{fig:temperature_tdi} compares results at $\mathit{Kn} = 0.001$ using $\lambda = 5$ and $\lambda = 500$. When $\lambda = 5$, the velocity nodes are overly concentrated, resulting in inaccurate predictions of temperature and flow fields. In contrast, $\lambda = 500$ produces results that closely match the analytical solution.

\begin{figure}[H]
	\centering
	\begin{subfigure}[b]{0.45\textwidth}
		\includegraphics[width=\textwidth]{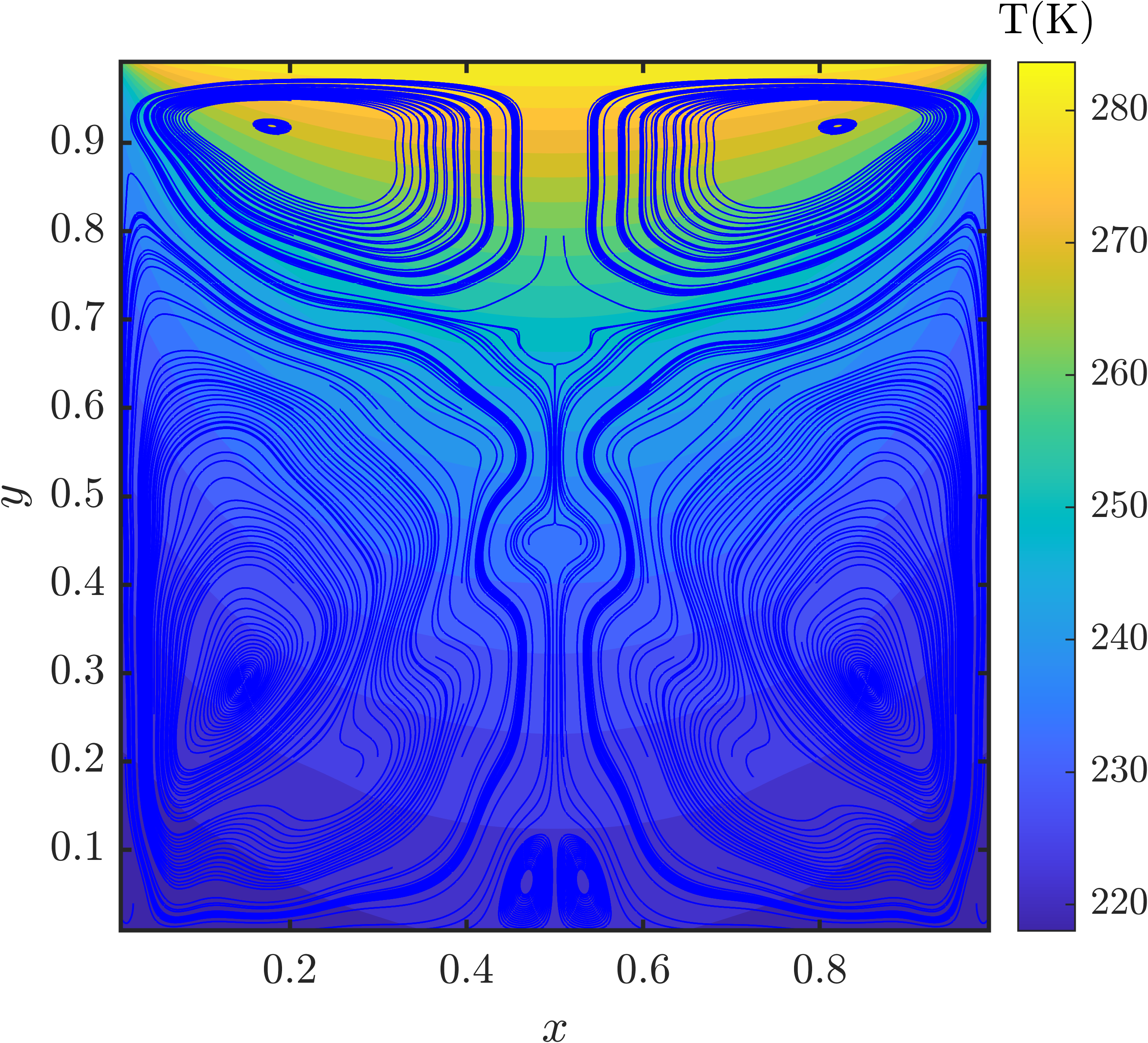}
		\caption{}
		\label{fig:NC}
	\end{subfigure}
	\quad
	\begin{subfigure}[b]{0.45\textwidth}
		\includegraphics[width=\textwidth]{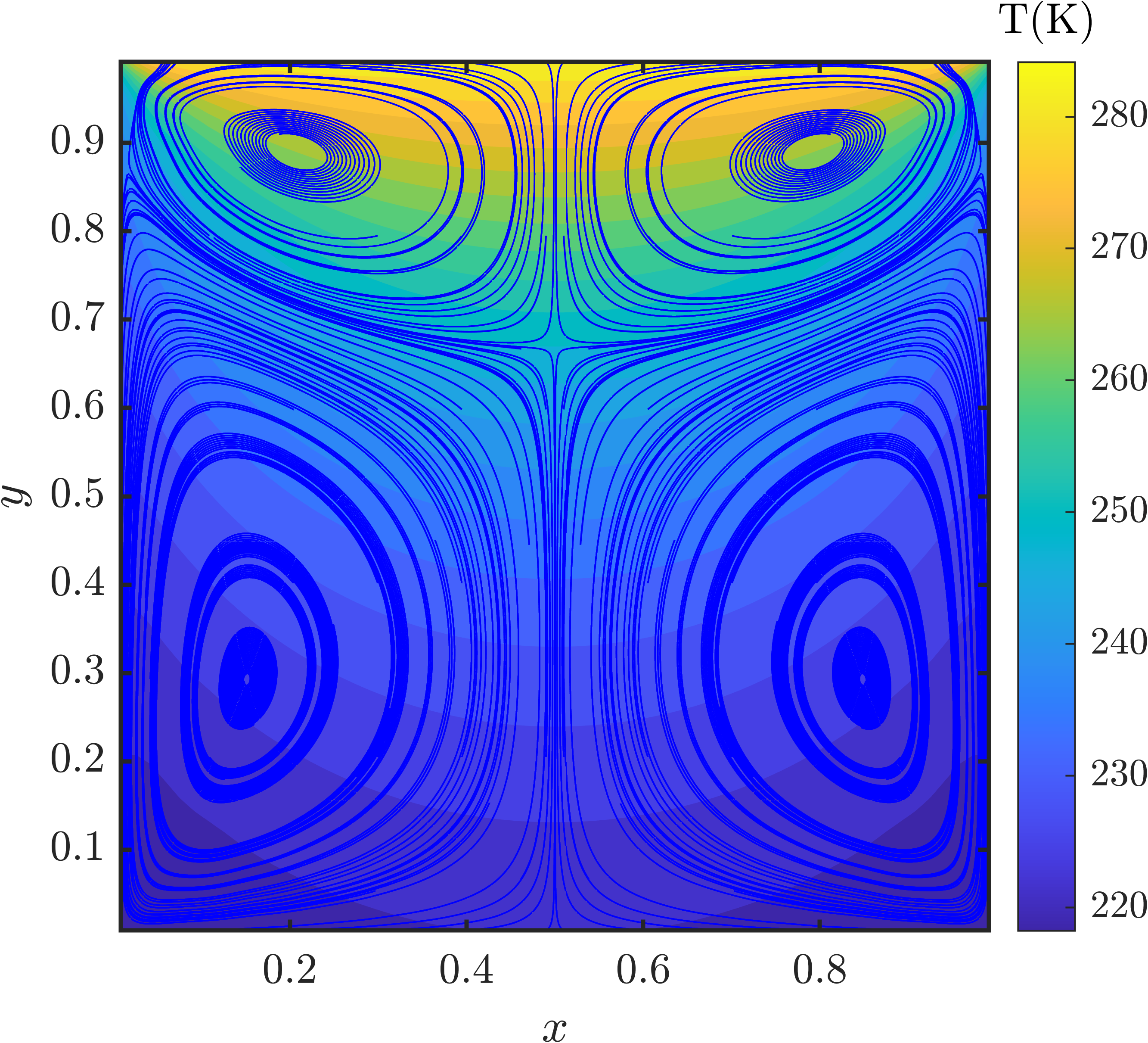}
		\caption{}
		\label{fig:ATGJ}
	\end{subfigure}
	\caption{Temperature field and streamlines of the TDI cavity flow at Kn = 10: 
		(a) NC rule with 201 × 201 nodes; 
		(b) ATGJ rule with 8 × 90 nodes.}
	\label{fig:TDI_Kn10}
\end{figure}

Fig.~\ref{fig:TDI_Kn10} presents the temperature fields and streamlines obtained by the Newton--Cotes (NC) quadrature ($201 \times 201$) and the present ATGJ method ($8 \times 90$). Despite the finer velocity grid, the NC solution exhibits unphysical vortices and lacks smoothness. The ATGJ method, with significantly fewer velocity points, achieves a smoother and more physical flow field, effectively eliminating the spurious structures observed in the NC results.

\begin{figure}[H]
	\centering
	\begin{subfigure}[b]{0.45\textwidth}
		\includegraphics[width=\textwidth]{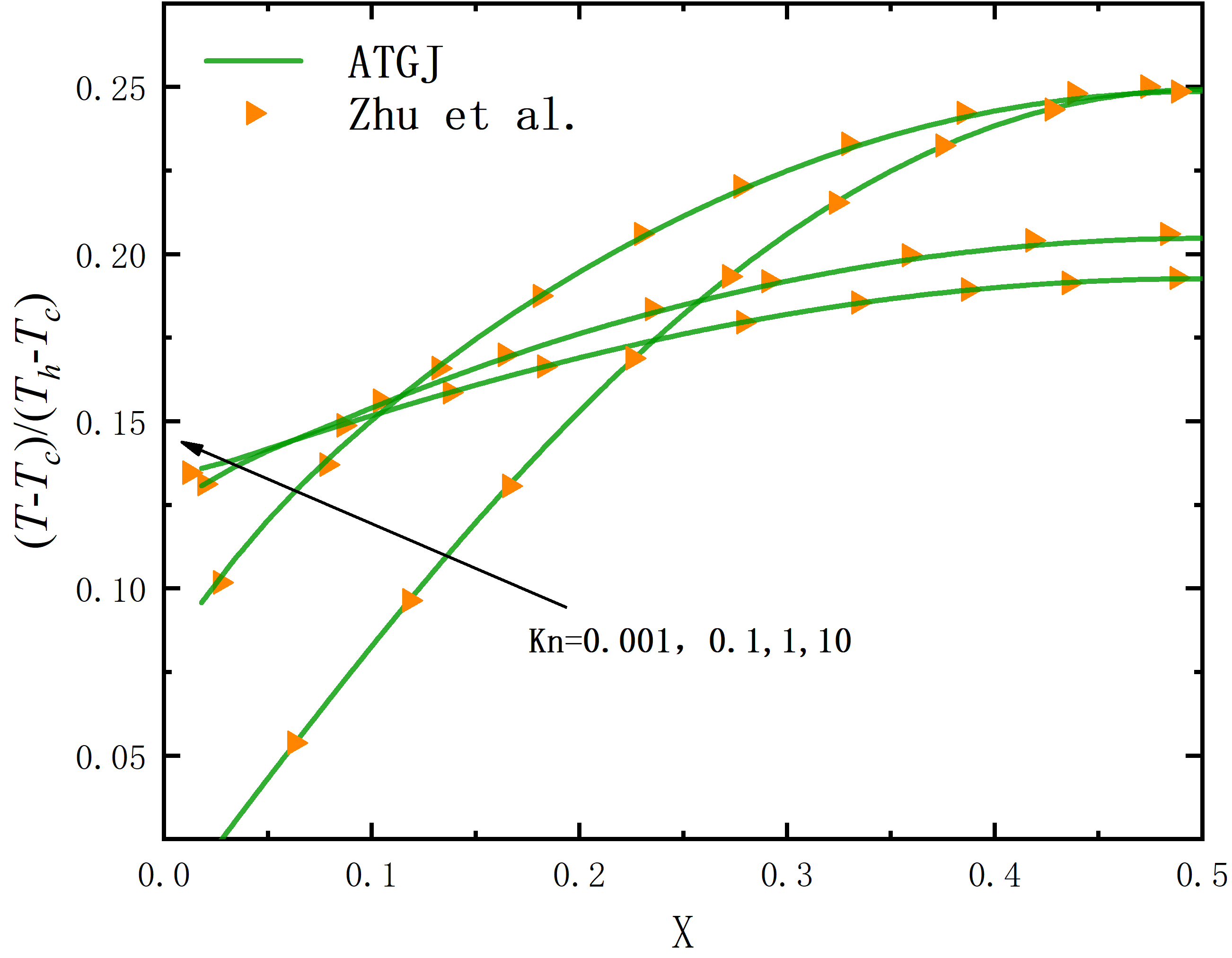}
		\caption{}
		\label{fig:o1}
	\end{subfigure}
	\hfill
	\begin{subfigure}[b]{0.45\textwidth}
		\includegraphics[width=\textwidth]{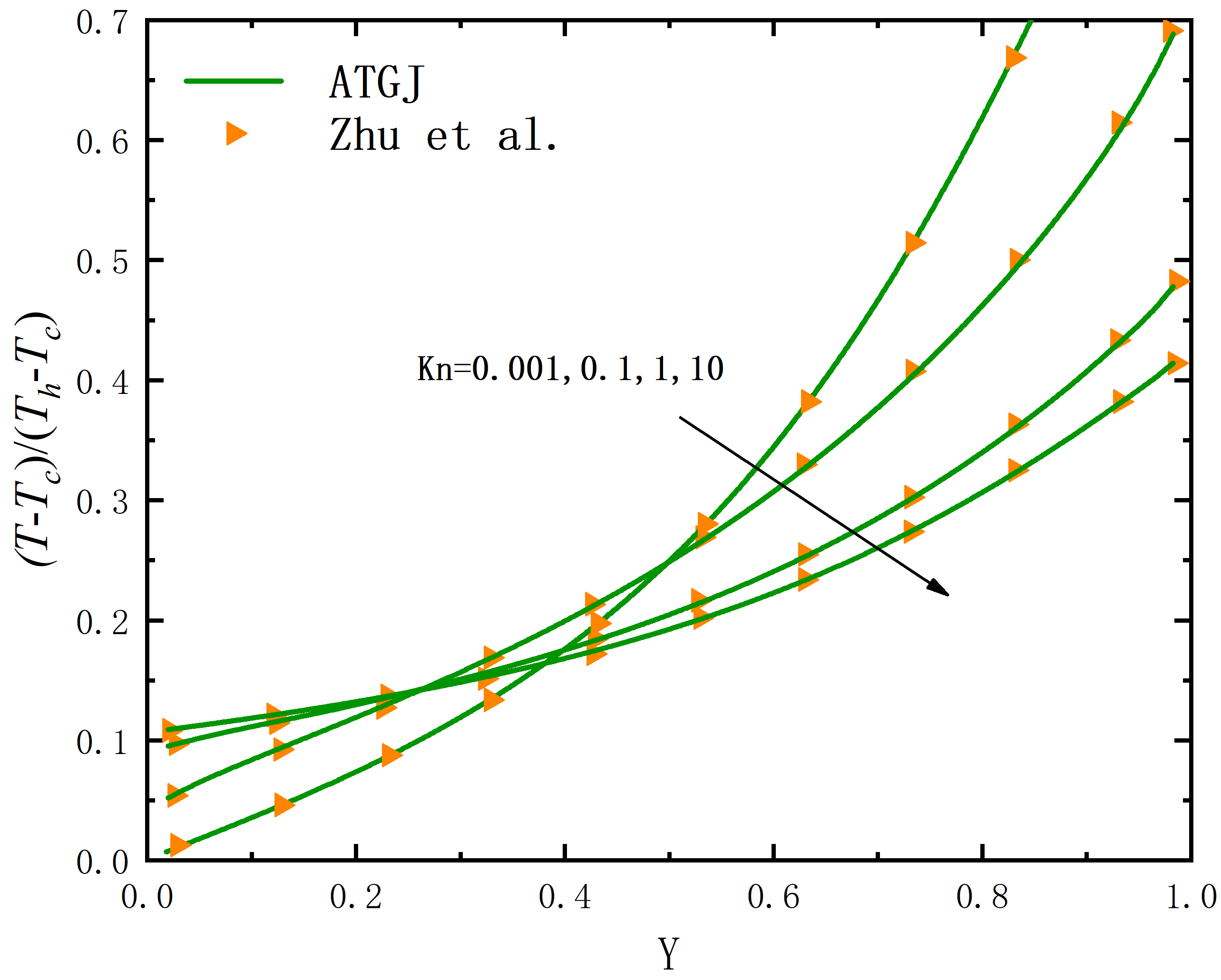}
		\caption{}
		\label{fig:o2}
	\end{subfigure}
	
	\vspace{2mm} % 行间距
	
	\begin{subfigure}[b]{0.45\textwidth}
		\includegraphics[width=\textwidth]{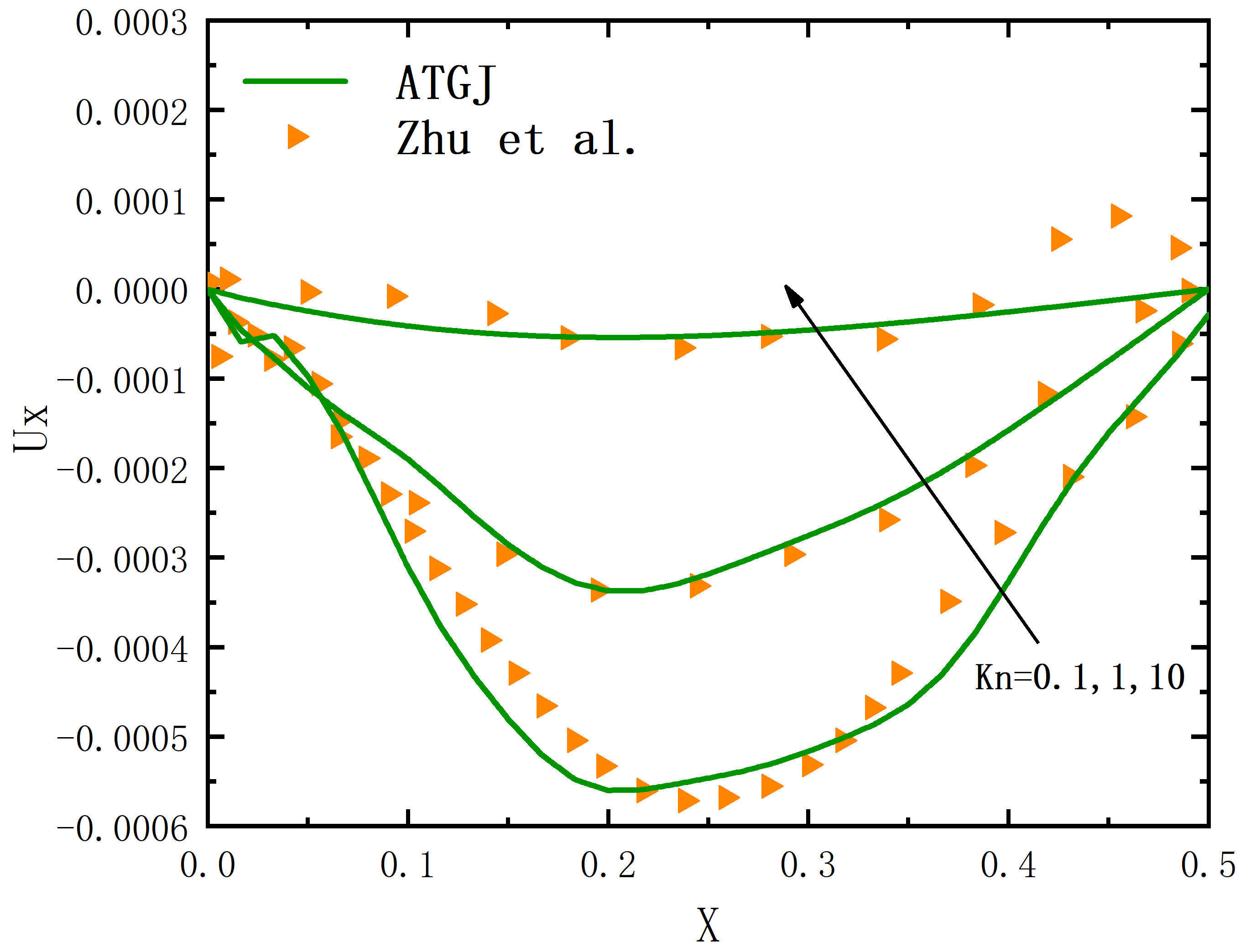}
		\caption{}
		\label{fig:o3}
	\end{subfigure}
	\hfill
	\begin{subfigure}[b]{0.45\textwidth}
		\includegraphics[width=\textwidth]{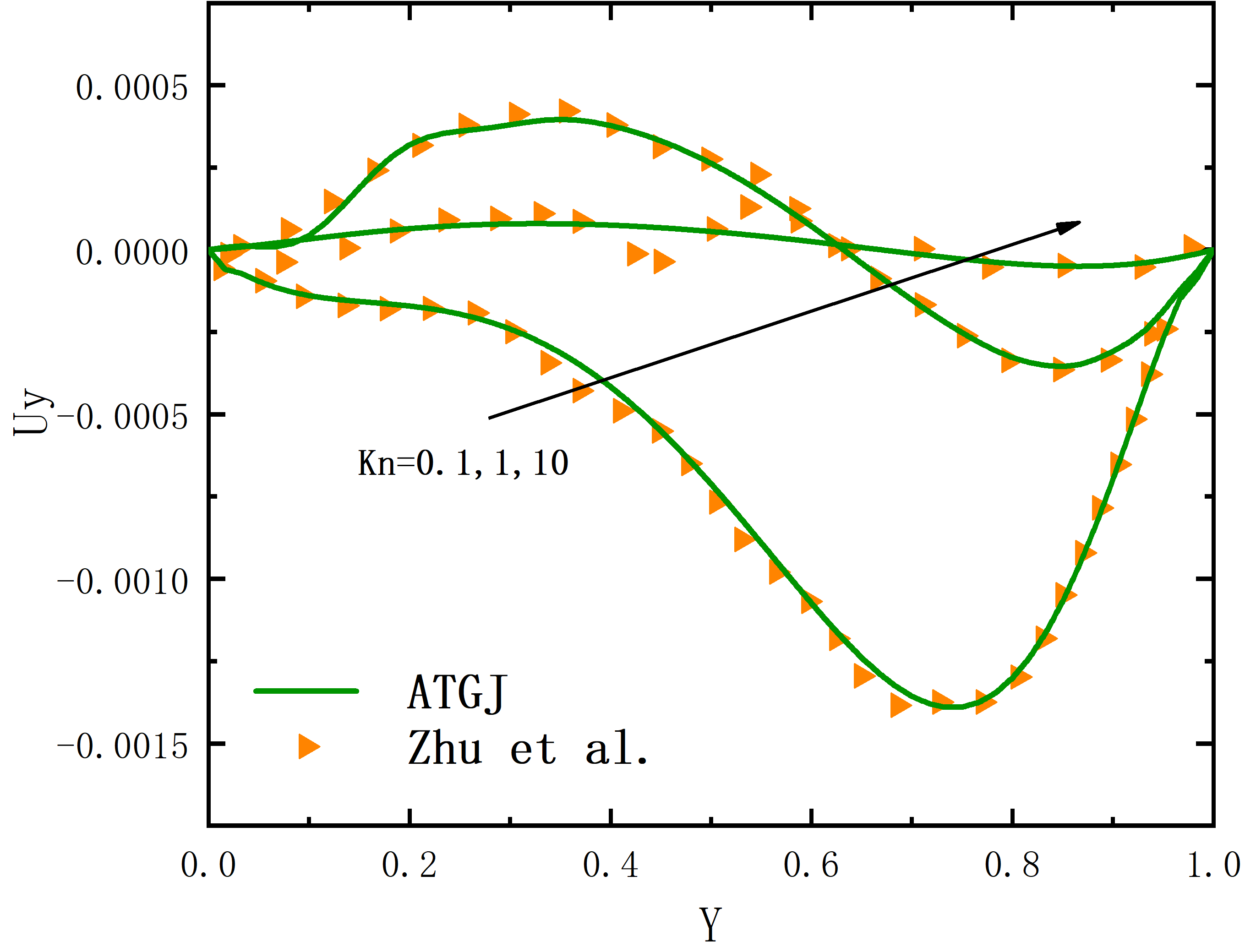}
		\caption{}
		\label{fig:o4}
	\end{subfigure}
	
	\caption{Temperature and velocity profiles along the horizontal and vertical centerlines for TDI cavity flow.}
	\label{fig:tdi_profiles}
\end{figure}

Fig.~\ref{fig:tdi_profiles} further compares temperature and velocity profiles along the cavity centerlines. Temperature predictions using ATGJ show excellent agreement with those reported by Zhu \textit{et al.} However, due to the small magnitude of velocities in this case, quantitative agreement in velocity profiles is more difficult to achieve. Notably, at $\mathit{Kn} = 10$, Zhu’s results display sharp peaks in both horizontal and vertical velocity at the cavity center, coinciding with the formation of artificial vortices in Fig.~4(a). The ATGJ results avoid these artifacts and yield smoother, more physically consistent profiles, demonstrating superior robustness and efficiency over conventional quadrature methods such as NC.
\subsection{Supersonic Flow past a Square Cylinder}
\label{subsec2}
This case considers a supersonic flow past a square cylinder at Mach number $Ma = 5$ and Knudsen number $Kn = 0.1$. The freestream conditions are prescribed as temperature $T_{\infty} = 1.0$, density $\rho_{\infty} = 1.0$, and velocity $u_{\infty} = 4.56$. The cylinder wall is maintained at a constant temperature $T_W = 1.0$, and its side length is $D = 1.0$. The computational geometry is shown in Fig.~2(b). Following the setup in Ref.~\citet{ref19}, the physical domain is discretized using 33,300 control volumes.

In velocity space, the ATGJ quadrature is employed with a $20 \times 60$ nodal configuration. For high-Mach-number flows, the discrete velocity domain must accommodate the freestream bulk velocity as well as thermal fluctuations, typically requiring coverage up to $u_{\infty} + 4\sqrt{RT}$. Accordingly, the parameters are set to $\alpha = 20$ and $\lambda = \frac{2\alpha}{\pi} + 20$, resulting in a velocity domain approximately bounded by a circle of radius 11. This ensures accurate resolution of the shifted and broadened distribution function associated with strong compressibility effects.

\begin{figure}[h]
	\centering
	\begin{subfigure}[b]{0.3\textwidth}
		\includegraphics[width=\textwidth]{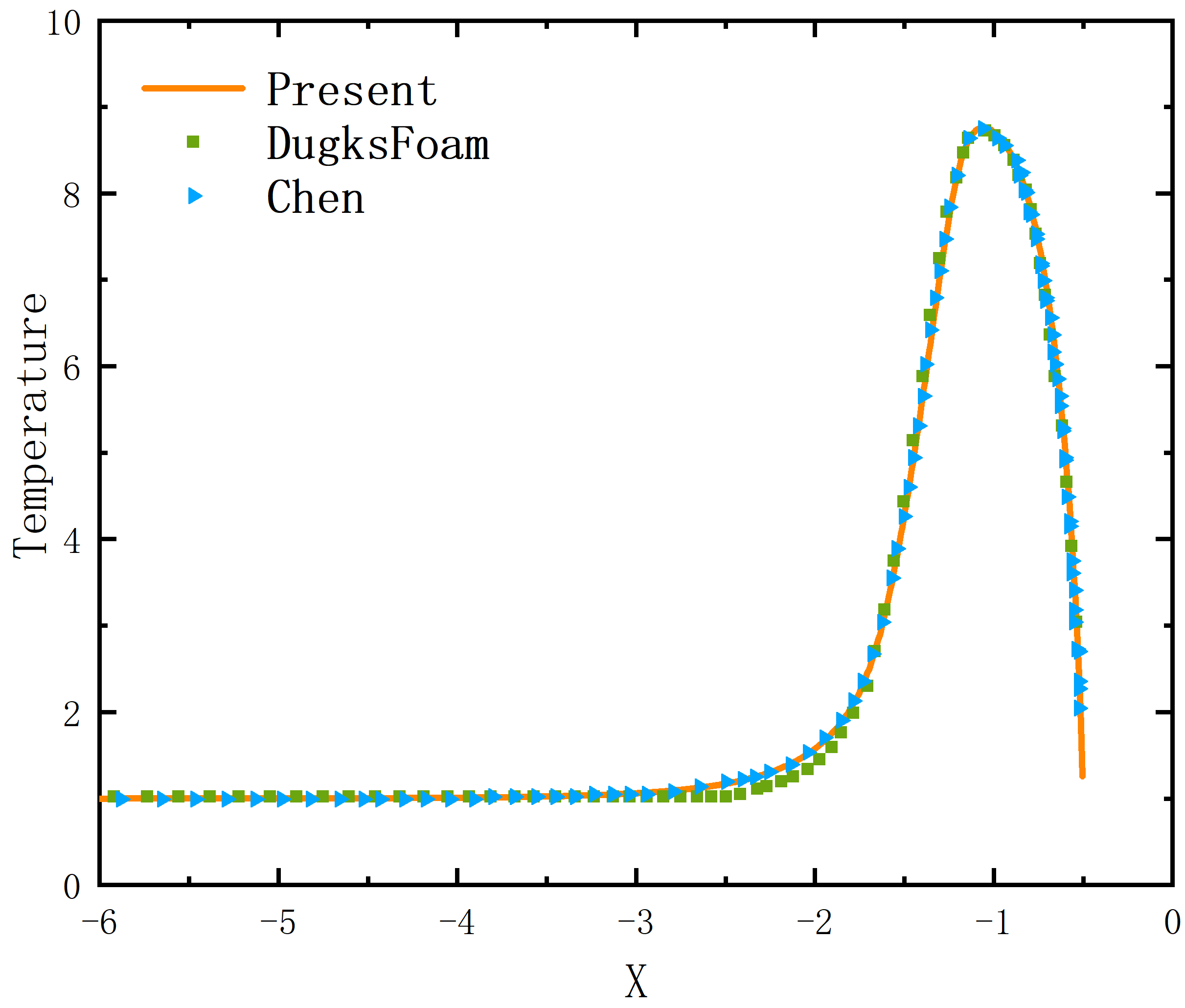}
		\caption{}
		\label{fig:s1}
	\end{subfigure}
	\hfill
	\begin{subfigure}[b]{0.3\textwidth}
		\includegraphics[width=\textwidth]{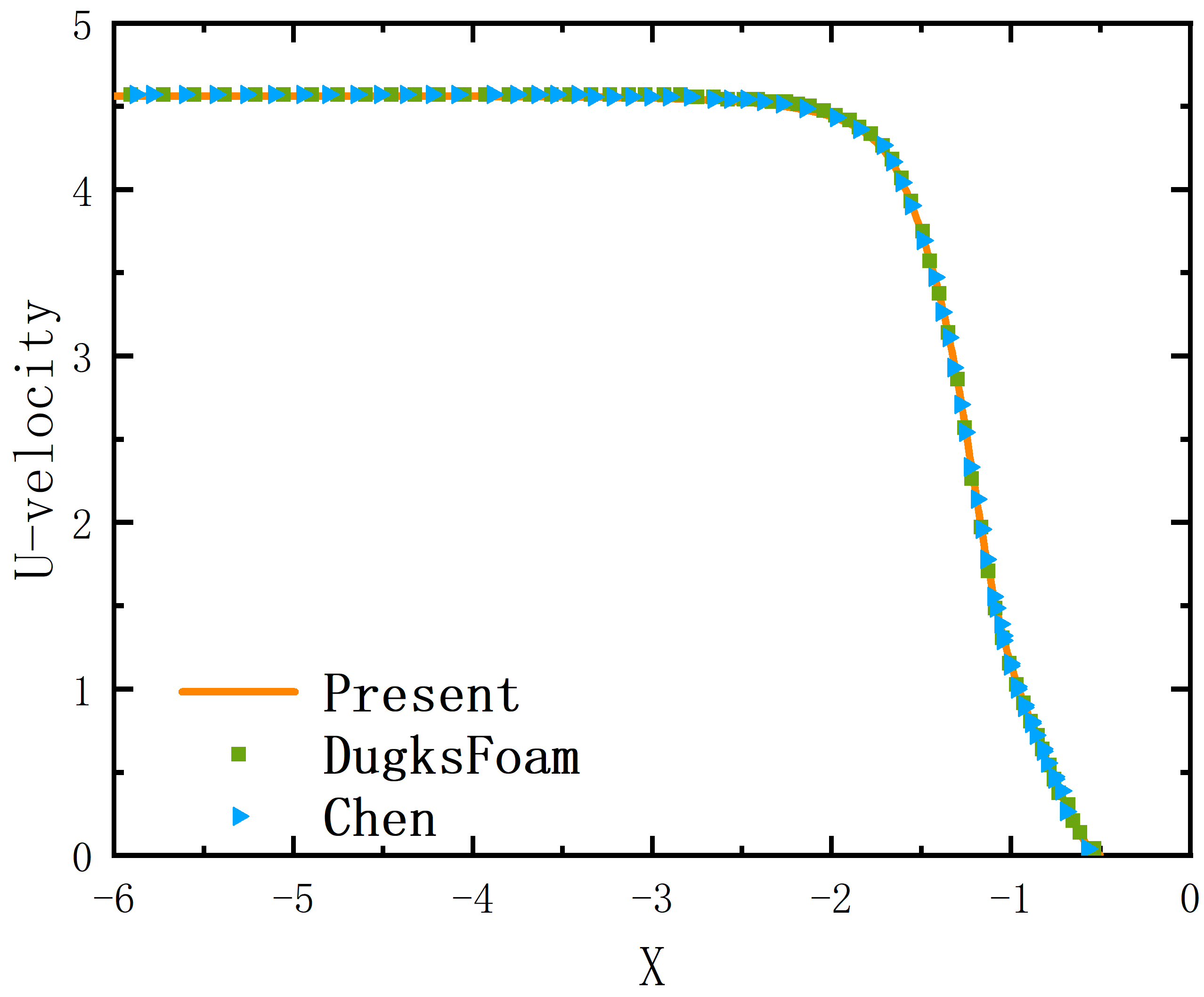}
		\caption{}
		\label{fig:s2}
	\end{subfigure}
	\hfill
	\begin{subfigure}[b]{0.3\textwidth}
		\includegraphics[width=\textwidth]{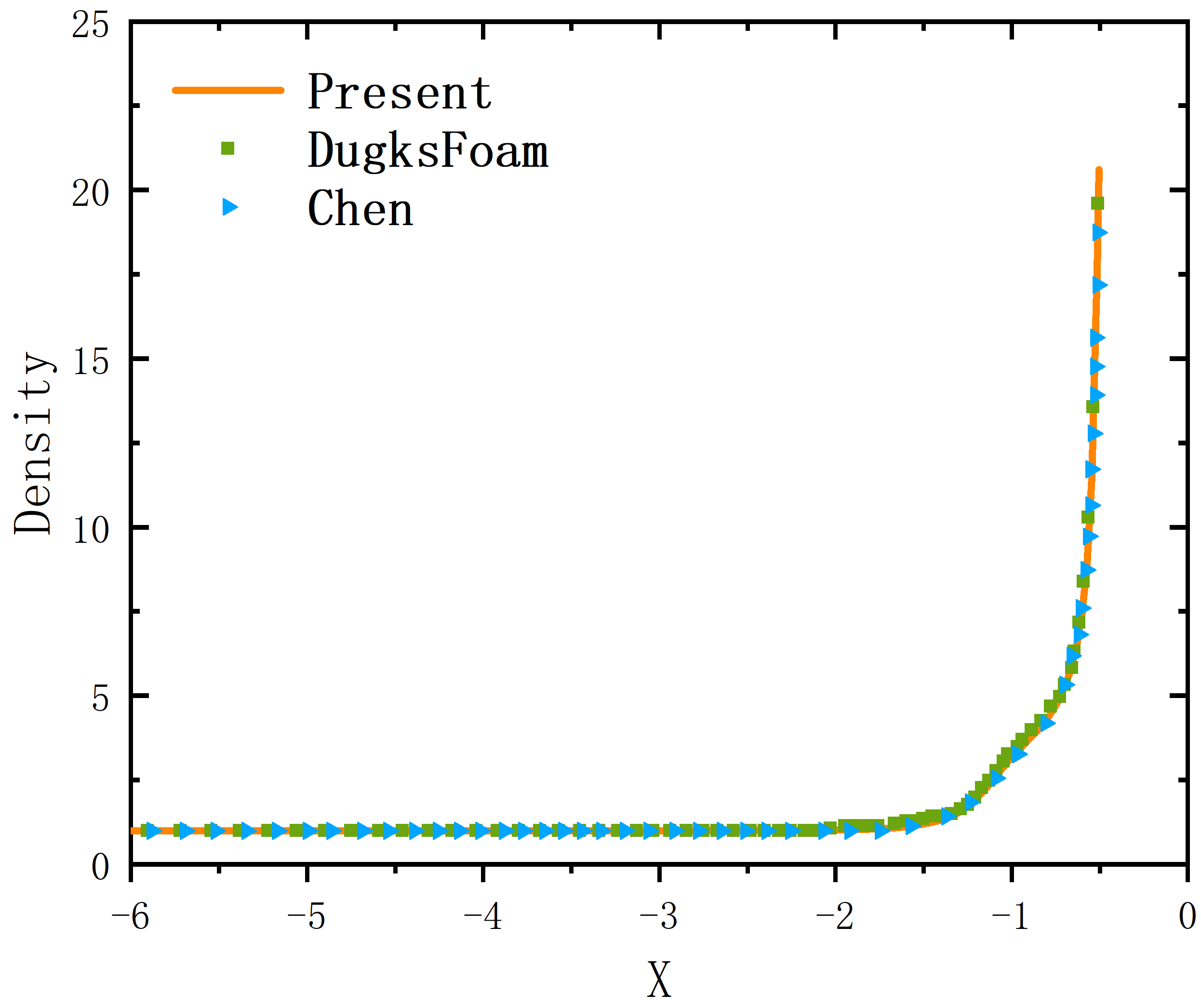}
		\caption{}
		\label{fig:s3}
	\end{subfigure}
	
	\caption{Distributions of temperature, velocity, and density flux along the center-symmetric line upstream of the stagnation point at $Ma = 5$.}
	\label{fig:ma5_distributions}
\end{figure}

Fig.~\ref{fig:ma5_distributions} shows the profiles of temperature, velocity, and density flux along the upstream centerline approaching the stagnation point. Results from the present ATGJ method are compared with those reported in Refs.~\citet{ref19,ref20}. Ref~\citet{ref19} utilizes the NC quadrature with $101 \times 101$ uniformly distributed velocity points. Despite using only 11\% of the discrete velocities, the ATGJ results exhibit excellent agreement with the reference data across all flow variables.

These results highlight the high accuracy and computational efficiency of the proposed quadrature in capturing shock structures, strong gradients, and rarefaction effects in supersonic and rarefied regimes. The ATGJ scheme proves to be a robust and cost-effective alternative to conventional uniform-grid quadrature methods in hypersonic gas dynamics simulations.
\section{Conclusion}
\label{sec5}
A novel Gaussian quadrature rule is developed in this work for velocity space discretization in multiscale Boltzmann solvers by introducing a parameterized weight function. The proposed scheme enables flexible control over the distribution of quadrature nodes and weights, thereby enhancing its ability to accurately capture non-equilibrium velocity distributions across a wide range of flow regimes. Numerical results demonstrate that the method offers improved accuracy, stability, and computational efficiency compared to conventional Newton--Cotes quadrature.

In particular, the scheme shows excellent adaptability in both continuum and rarefied regimes, making it well-suited for multiscale flow simulations. Owing to its analytical construction and tunable parameters, the method exhibits good generality and can be extended to multi-dimensional velocity spaces. Future work will explore automated parameter selection and three-dimensional implementations for more complex flow applications.

\section{Acknowledgements} G.Q.W. acknowledges support from the Pioneer R\&D Program of Zhejiang Province-China (Grant No.~2023C03016).

\appendix
\section*{Appendix: Core Function Implementation}
\label{app:dvATGJ}

For the two-dimensional Gauss--Jacobi quadrature proposed in this work, 
the discrete velocities and weights are computed using Eq.~\eqref{eq:discrete_nodes}, 
where $\theta_j$ and $w_{\theta_j}$ are given by Eq.~\eqref{eq:theta_nodes}. 
The radial quadrature nodes and the corresponding weights are available at 
\url{https://github.com/WangLu521/2DATGJ.git}.
 
Below, we provide the core Python function used to generate the quadrature points and weights.

\begin{lstlisting}[caption={Core function \texttt{dvATGJ}}, label={lst:dvATGJ}]
import numpy as np
from scipy.linalg import eigh
from scipy.special import gammaln

def rootsWeights(n, alpha, lamda, T0):
  if n <= 0 or alpha <= 0 or lamda <= 0:
    raise ValueError("n, alpha, lamda must be > 0.")

  beta = 0
  a, b = jacobi_recurrence(n, alpha, beta)
  J = np.diag(a) + np.diag(b[1:], 1) + np.diag(b[1:], -1)
  eigenvalues, eigenvectors = eigh(J)

  roots = 0.5 * (eigenvalues + 1)
  Rr = np.sqrt(lamda * T0 * np.tan(np.pi/2 * roots))

  log_ratio = gammaln(alpha + 1) - gammaln(alpha + 2)
  weights = np.exp(log_ratio) * (eigenvectors[0, :]**2)

  return Rr, weights
\end{lstlisting}

**Parameter description:**  

\begin{itemize}
	\item \texttt{n}: Number of quadrature points along the radial direction.
	\item \texttt{alpha}: Shape parameter of the Jacobi polynomial.
	\item \texttt{lambda}: Scaling factor used in the transformation.
	\item \texttt{T0}: Reference temperature parameter.
	\item Returns: Radial quadrature points \texttt{Rr} and corresponding \texttt{weights}.
\end{itemize}

\begin{lstlisting}
def jacobi_recurrence(n, alpha, beta):
  a = np.zeros(n)
  b = np.zeros(n)
  for k in range(n):
    if k == 0:
      a[k] = (beta - alpha) / (alpha + beta + 2)
    else:
      a[k] = (beta**2 - alpha**2) / ((2 * k + alpha + beta) * (2 * k + alpha + beta + 2))
    if k > 0:
      num = 4 * k * (k + alpha) * (k + beta) * (k + alpha + beta)
      den = (2 * k + alpha + beta)**2 * (2 * k + alpha + beta + 1) * (2 * k + alpha + beta - 1)
      b[k] = np.sqrt(num / den)
    return a, b
\end{lstlisting}

\end{document}